\documentclass{amsart}

\usepackage{epstopdf}
\usepackage{graphicx}
\usepackage{psfrag}
\usepackage{color,xcolor}

\theoremstyle{cupthm}
\newtheorem{theorem}{Theorem}[section]

\newtheorem{lemma}[theorem]{Lemma}
\theoremstyle{cupdefn}

\theoremstyle{cuprem}
\newtheorem{remark}[theorem]{Remark}
\numberwithin{equation}{section}

\newcommand{\R}{{\mathbb R}}
\newcommand{\Z}{{\mathbb Z}}
\newcommand{\B}{{\mathcal B}}
\newcommand{\E}{{\mathbb E}}
\def\H{\mathcal{H}}

\def\vep{\varepsilon}
\def\gtr{{\mathcal{G} T_r}}
\def\gtrss{{\mathcal{G} T_r^*}}
\def\g{\mathcal{G}}
\def\hath{\widehat{\mathcal{H}}}

\newcommand{\leb}{{\rm Leb\,}}
\newcommand{\card}{{\rm card\,}}

\begin{document}
\title{Local level sets of  the Takagi-van der Waerden function}
\author{Lai Jiang}
\address{Hangzhou Institute for Advanced Study, UCAS}
\email{jianglai@ucas.ac.cn}

\author{Ting-Ting Ying}
\address{Zhejiang University}
\email{3220102483@zju.edu.cn}

\author{Yi-Yang Zhang}
\address{Zhejiang University}
\email{3220102591@zju.edu.cn}




\begin{abstract}

In this paper, we investigate the Takagi-van der Waerden function,
$$
	T_r(x) = \sum_{n=0}^{\infty} \frac{\phi(r^n x)}{r^n} ,\quad x\in [0,1], \quad r \in \mathbb{Z}^+,
$$
where $\phi(x)={\rm dist}(x,\mathbb{Z})$ represents the distance from $x$ to the nearest integer. 

Lagarias and Maddock [Level sets of the Takagi function: local level sets, \emph{Monatsh. Math.}, {\bf  166} (2012), No. 2, 201--238] introduced the notion of local level sets for the classical Takagi function $T_2$. They proved that if $y$ is a random variable uniformly distributed over the range of $T_2$, then the expected number of local level sets contained in the level set $L_2(y)$ equals $3/2$.
We extend the study by defining an analogous concept of local level sets for all even integers $r$. Then we prove that, for every even integer $r\geq 2$, if $y$ is a random variable uniformly distributed, then the expected number of local level sets contained in the level set $L_r(y)$ equals $1 + 1/r$.

\end{abstract}


\maketitle

\section{Introduction}

In 1903, Takagi constructed a continuous but nowhere differentiable function using the following formula:
$$
T(x) = \sum_{n=0}^{\infty} \frac{\phi(2^n x)}{2^n}, \quad x \in [0,1],
$$
where $\phi(x) = {\rm dist}(x, \Z)$ is the distance from $x$ to the nearest integer. This function is regarded as the classical Takagi function. A well-known generalization is called the Takagi-van der Waerden function, which is defined as
$$
T_r(x) = \sum_{n=0}^{\infty} \frac{\phi(r^n x)}{r^n}, \quad x \in [0,1], \quad \forall r \geq 2.
$$

It is straightforward to conclude that $T_r$ is also continuous and nowhere differentiable. 
There are many works on the image of $T_r$; see \cite{AK11} and the references therein.
For instance, Baba \cite{B84} studied the maximum value of $T_r$. Concerning dimensions, both the Hausdorff dimension and box dimension of the graph of $T_r$ are one. Additionally, the first author \cite{J25} demonstrated that its Assouad dimension equals one.

When a line intersects a set, the resulting intersection is referred to as a slice. Many important problems fall into this category, including Furstenberg's conjecture \cite{W19}. Anttila, Bárány, and Käenmäki \cite{ABK24} explored the connection between slices and dimensions. For any $r \geq 2$ and $y \in \R$, we define
$$
L_r(y) := \big\{ x \in [0,1] : T_r(x) = y \big\}.
$$
We call $L_r(y)$ the \emph{level set} of $T_r$, which is a special class of slices.

A natural question is whether the level set $L_r(y)$ is finite or infinite. There are many interesting results concerning the level sets of the Takagi function \cite{A12,A13,B08}, particularly in the case where $r = 2$. Buczolich \cite{B08} showed that $L_2(y)$ is finite for almost every $y$.

Lagarias and Maddock \cite{LM11,LM12} deeply studied the properties of level sets. They classified the points in level sets by their expansions.
For any $x \in (0,1)$, its binary expansion is given by
$$
x = 0.\vep_1 \vep_2 \cdots \vep_k \cdots = \sum_{k=1}^\infty \frac{\vep_k}{r^k}, \quad \vep_k \in \{0, 1\}, \quad \forall k \in \Z^+.
$$
Let $\card(\cdot)$ be the cardinality of a set.
For any $n \in \Z^+$, define
$$
D_{2,n}(x) := \card \big\{ 1 \leq j \leq n : \vep_j = 0 \big\} - \card \big\{ 1 \leq j \leq n : \vep_j = 1 \big\}.
$$
They found that if $x,x' \in (0,1)$ satisfy $|D_{2,j}(x) |=|D_{2,j}(x')|$ for all $j \in \Z^+$, then $T_2(x)=T_2(x')$. From this property, they defined an equivalence relation on $(0,1)$ by
$$
x \sim_2 x' \iff |D_{2,j}(x)| = |D_{2,j}(x')|\quad \text{for all }j \in \Z^+.
$$
Based on the equivalence relation, Lagarias and Maddock \cite{LM12} introduced the concept of \emph{local level sets} for the classical Takagi function $T_2$, which is
$$L_{2,x}^{loc} = \{x' : x' \sim_2 x\}.$$
They also proved the following interesting result.

\begin{theorem}[\cite{LM12}]\label{thm-0}

Let $y$ be a random variable uniformly distributed on the range of $T_2$. Then the expected number of local level sets contained in $L_r(y)$ is given by
$$
	\E \big[N_2^{loc}(y) \big] = \frac{3}{2}\int_0^{\frac{2}{3}} N_2^{loc}(y) \, dy = \frac{3}{2}.
$$
Here $N_2^{loc}(y)$ is the number of local level sets contained in $L_2(y)$. 
\end{theorem}

Moreover, they pointed out that each local level set is either finite or a Cantor set. We remark that the local level sets are much easier to analyze and provide methods for dissecting level sets and understanding their structure.
It is also worth mentioning that local level sets have a strong and direct connection to random walks \cite{A12,LM12}.

In this paper, we generalize Theorem~\ref{thm-0} to all even integers $r$. First we define an equivalence relation $\sim _r$ as follows.
Given $x \in (0,1)$, its $r$-based expansion is given by
$$
x = 0.\vep_1 \vep_2 \cdots \vep_k \cdots = \sum_{k=1}^\infty \frac{\vep_k}{r^k}, \quad \vep_k \in \{0, 1, \ldots, r-1\}, \quad \forall k \in \Z^+.
$$
For any $n \in \Z^+$, we define
$$
D_{r,n}(x) := \card \big\{ 1 \leq j \leq n : \vep_j < r/2 \big\} - \card \big\{ 1 \leq j \leq n : \vep_j \geq r/2 \big\} . 
$$

Assume that $x = 0.\vep_1 \vep_2 \ldots$ and $x' = 0.\eta_1 \eta_2 \ldots$ are two different elements in $(0,1)$. We say $x \sim_r x'$ if the following conditions are satisfied for each $j \in \Z^+$:
\begin{itemize}
    \item[(C1)] $|D_{r,j}(x)| = |D_{r,j}(x')|$,
    \item[(C2)]  either $\vep_j = \eta_j$ or $\vep_j + \eta_j = r - 1$ if $|D_{r,j}(x)|>0$.
\end{itemize}

We will see that $T_r(x)=T_r(x')$ for all $x\sim_rx'$ in Lemma \ref{lem:sim}, which motivates our definition of local level sets:
$$
	L_{r,x}^{loc} = \{x' : x' \sim_r x\}.
$$
Notice that in the case $r=2$, we always have either $\vep_j=\eta_j$ or $\vep_j+\eta_j=1$.
Thus, our definition coincides with the original one for $r=2$. 

We remark that condition (C2) is necessary since condition (C1) holds for $x=0.4$ and $x'=0.9$ in the case $r=10$, but $T_{10}(x) >T_{10}(x')$.

Our main result extends Theorem~\ref{thm-0} to all even integers $r$. 

\begin{theorem}\label{thm-1}
Let $r \geq 2$ be an even integer and $y$ be a random variable uniformly distributed on the range of $T_r$. Then the expected number of local level sets contained in $L_r(y)$ is given by 
$$
	\E \big[N_r^{loc}(y) \big] = \frac{2r^2-2}{r^2}\int_0^{\frac{r^2}{2r^2-2}} N_r^{loc}(y) \, dy = 1 + \frac{1}{r}.
$$
Here $N_r^{loc}(y) $ is the number of local level sets contained in $L_r(y)$.
\end{theorem}

We also prove that the expected number of elements in $L_r(y)$ is infinite.

The paper is organized as follows. 
In Section 2, we present some basic properties of Takagi-van der Waerden function and local level sets. 
In Section 3, we recall the definition of humps and obtain the cardinality of humps. 
In Section 4, we derive the expected number of local level sets and elements contained within level sets.

\section{Basic properties}

For any $n \in \Z^+$, we define the partial sum sequences of $T_r$ as
\begin{equation}\label{eq:partial-def}
	S_{r,n}(x) := \sum_{k=0}^{n-1} \frac{\phi(r^k x)}{r^k}, \quad x \in [0,1].
\end{equation}
It is easy to conclude that $S_{r,n}$ uniformly converges to $T_r$.

\begin{lemma}\label{lem:graph-simi}
For any integers $n \in \Z^+$ and $0 \leq k \leq r^n-1$, we have 
\begin{equation}\label{eq:simi-value}
T_r \big(\frac{k}{r^n}+t \big) = T_r\big(\frac{k}{r^n}\big) + \frac{1}{r^n} T_r(r^n t ) + D_{r,n}\big(\frac{k}{r^n}\big) \cdot t, 
\quad \forall t \in [ 0, r^{-n}].
\end{equation}
\end{lemma}
\begin{proof}
Fix integers $n$ and $k$. 
Write $x=k/r^n$ and assume that the $r$-based expansion of $x$ is given by $x = 0.\vep_1 \cdots \vep_n $.
Then for any $t \in [0,r^{-n}]$, we have
\begin{equation}\label{eq:value-i-1}
	T_r(x+t) -S_{r,n}(x+t)
	=\sum_{i=n}^\infty \frac{\phi(r^i (x+t))}{r^i} 
	=\frac{1}{r^n} \sum_{i=n}^\infty \frac{\phi(r^{i-n} \cdot r^n t )}{r^{i-n}} 
	= \frac{1}{r^n} T_r( r^n t) .
\end{equation}

Fix integer $1 \leq i \leq n$, write $a_i=1$ if $\vep_{i} <r/2 $, $a_i=-1$ if $\vep_i \geq r/2 $.
Since $ \phi(z) $ is linear on both $[0,1/2]$ and $[1/2,1]$, 
then $ \phi(r^{i-1} z) $ is linear on $ [ x, x+r^{-n}] $,
$$
	\frac{\phi(r^{i-1}(x+t))- \phi(r^{i-1} x)}{r^{i-1}} = \frac{a_i (r^{i-1}(x+t)-r^{i-1} x)}{r^{i-1}}= a_i t.
$$
Summing over $i$ from $1$ to $n$, we have $D_{r,n}(x)=\sum_{i=1}^{n} a_i$ and
\begin{equation}\label{eq:value-i-2}
	S_{r,n}(x+t)
	=\sum_{i=0}^{n-1} \Big( \frac{\phi(r^{i-1}x)}{r^{i-1}} + a_i t \Big)
	=T_r (x) + D_{r,n}(x) \cdot t.
\end{equation}
Combining \eqref{eq:value-i-1} and \eqref{eq:value-i-2} completes the proof.
\end{proof}

For any $x \in [0,1)$, $x$ is called an $r$-based rational if $r^n x \in \Z$ for some $n \in \Z^+$.
Moreover, an $r$-based rational $x$ is called \emph{balanced} if there exists $n \in \Z^+$ such that $r^{n}x \in \Z$ and $D_{r,n}(x) = 0$. 
If $r^p x \in \Z$, then $D_{r,p+j}(x)=D_{r,p}(x)+j$ for all $j \in \Z^+$.
This implies that $n$ is unique and even. 
We write $I(x)=[x,x +r^{-n}]$.

Write $\B^r$ for the set of all $r$-based balanced rationals and
$$ \B_1^r= \{ x \in \B^r: \mbox{ there exists unique } j \in \Z^+ \mbox{ such that } D_{r,j}(x) =0 \} .$$
For any $ x =0. \vep_1 \vep_2 \cdots \in \B^r $, let $n= \min\{ j \in \Z^+: D_{r,j}(x) =0 \}$ and $x_0=0. \vep_1 \cdots \vep_n$.
It is clear that $x_0 \in \B_1^r$ and $I(x) \subset I(x_0)$.
As a result, 
$$
	\bigcup_{x_0 \in \B^r} I(x_0) =  \bigcup_{x_0 \in \B_1^r} I(x_0).
$$
In \cite{AK11}, $\B_1^2$ is called the set of binary balanced points of first generation.

Write $\gtr:=\{(x,T_r(x)):x\in [0,1]\}$ and $\gtrss =\gtr |_{X_r^*}$, where
$$
	X_r^* : = [0,1] \setminus \bigcup_{x_0 \in \B^r} I(x_0) = [0,1] \setminus \bigcup_{x_0 \in \B_1^r} I(x_0).
$$
Then for any $x \in (0,1)$, $x \in X_r^*$ if and only if $D_{r,j}(x) \neq 0$ for all $j \in \Z^+$.

\begin{remark}

Notice that $\g T_r$ is symmetric in the sense that $T_r(z) = T_r(1-z)$ for all $z \in [0,1]$. 
Assume that $x \in \B^r$ with $D_{r,n}(x)=0$ and $r^n x \in \Z$, then for any $t \in [0, 1/r^n]$,
\begin{equation}\label{eq:simi-eq}
T_r(x+t) = T_r(x) + \frac{ 1 }{r^n} T_r(r^n t) = T_r(x) + \frac{ 1 }{r^n} T_r(1-r^n t) = T_r (x + r^{-n} - t) .
\end{equation}
\end{remark}

\begin{lemma}\label{lem:tmp}
For any distinct $x',x'' \in X_r^* $, $T_r(x') = T_r(x'')$ if and only if $x'+x''=1$.
\end{lemma}

\begin{proof}
Define
$$
S_r^*(x) =
\begin{cases}
T_r(x_0), & \mbox{if } x \in I(x_0),\mbox{ where } x_0 \in \B_1^r, \\
T_r(x), &  \mbox{otherwise}.
\end{cases}
$$
By the definition of $S_r^*$, we know $T_r(x)=S_r^*(x)$ for all $x \in X_r^*$.

First, we are going to show that $S_r^*$ is increasing on $[0,1/2]$.
For any $m \in \Z^+$, write 
$$\Omega_m : =\{ x \in \B_1^r : |I(x)| \geq r^{-m}  \}. $$
Then for any $x \in \Omega_m$, we have $r^m x \in \Z$ and $ D_{r,m}(x) \geq 0 $.
Set
$$
	S_{r,m}^*(x) =
	\begin{cases}
		T_{r}(x_0), & \text{if } x \in I(x_0) \mbox{ for some } x_0 \in \Omega_m , \\
		S_{r,m}(x), & \mbox{otherwise},
	\end{cases}
$$
where $S_{r,m}$ is defined in \eqref{eq:partial-def}.

Fix an integer $0 \leq i < r^m/2$ and write $x_1=i/r^m$.
If $x_1 \in I(x_0) $ for some $x_0 \in \Omega_m$, we have $S_{r,m}^*(x)=T_r(x_0)$ for all $x \in [x_1,x_1+r^{-m}]$.
Otherwise, we have $D_{r,m}(x_1)>0$. By \eqref{eq:value-i-2}, $S_{r,m}$ is strictly increasing on $[x_1, x_1+r^{-m}]$.
Hence, $S_{r,m}^*$ is increasing on $[x_1,x_1+r^{-m}]$. 
Since $i$ is arbitrary and $S_{r,m}^*$ is continuous, $S_{r,m}^*$ is increasing on $[0,1/2]$.
Note that $S_{r,m}^*$ uniformly converges to $S_r^*$, so $S_r^*$ is increasing on $[0,1/2]$.

Pick $x_1,x_2 \in X_r^* \cap [0,1/2]$.
Then there exists $x_0 \in \B_r^*$ such that $x_1<x_0<x_2$. 
By the definition of $S_r^*$, $T_r(x)=S_r^*(x)$ holds for all $x \in X_r^*$.
Since $S_r^*$ is increasing, it follows that $S_r^*(x_1) \leq S_r^*(x_0) $.
From \eqref{eq:simi-value}, $T_r(x_0) < T_r(x_2)$.
Thus, $ T_r(x_1) \leq T_r(x_0) < T_r(x_2)$.
By combining with the symmetry of $\gtr$, the proof is complete.
\end{proof}

\begin{lemma}\label{lem:sim}
If $ x' \sim_r x$, then $T_r(x') = T_r(x)$. 
\end{lemma}
\begin{proof}
Without loss of generality, assume $x<x'$ and let
$x = 0.\vep_1  \cdots \vep_n \cdots$ and $x' = 0.\vep_1'  \cdots \vep_n' \cdots$, respectively. 
Write
$
	\Omega= \{ j \in \Z^+: D_{r,j}(x) = 0 \} .
$

If $\Omega$ is empty, we have $\vep_j =r-1-\vep_j'$ for all $j \in \Z^+$. This implies $x'=1-x$ and consequently $T_r(x)=T_r(x')$.

If $\Omega$ is nonempty, let $k=\min \Omega$.
Write $t= 0.\vep_1 \cdots \vep_k$ and $t'=0.\vep_1'  \cdots \vep_k'$, then $t \leq t'$. 
We will show that $T_r(t)=T_r(t')$.
From $x \sim_r x'$, we have either $\vep_1=\vep_1'$ or $\vep_1 +\vep_1'=r-1$.
In the case where $\vep_1=\vep_1'$, we have $D_{r,j}(x)=D_{r,j}(x')$ and $\vep_j=\vep_j'$ for all $1 \leq j <k$.
Moreover, $\vep_k '<r/2$ if and only if $\vep_k<r/2$.
By \eqref{eq:simi-eq},  
$$
	T_r(t)=T_r(t+r^{-n})= \cdots = T_r(t').
$$
If $\vep_1 \neq \vep_1'$, we set $t''=0.\vep_1'' \cdots \vep_k''$ where $\vep_j''=r-1-\vep_j'$ for all $1 \leq j \leq k$.
Clearly, $t'' \sim_r t' $ and
$$
	T_r(t')=T_r(t'+r^{-k})=T_r(1-r^{-k} -t')=T_r(t'').
$$
As a result, $t'' \sim_r t$ and $\vep_1 ''=\vep_1$. This implies $T_r(t)=T_r(t'')=T_r(t')$. 

Write $x_1 = 0.\vep_{k+1} \vep_{k+2}  \cdots$ and $x_1' = 0.\vep_{k+1}'   \vep_{k+2}' \cdots$, then $ x_1 \sim_r x_1'$.
By \eqref{eq:simi-eq},
$$
	T_r(x)-T_r(x')= \Big( T_r(t) + \frac{ T_r( x_1)}{r^k} \Big) - \Big( T_r(t') + \frac{ T_r( x_1')}{r^k} \Big)
	= \frac{ T_r( x_1)-T_r( x_1')}{r^k}.
$$
By induction, for any $p \in \Omega$ and $x_p = 0.\vep_{p+1} \vep_{p+2}  \cdots$ and $x_p' = 0.\vep_{p+1}'   \vep_{p+2}' \cdots$,
\begin{equation}\label{eq:tmp-2}
		T_r(x)-T_r(x')=  \frac{ T_r( x_p)-T_r( x_p')}{r^p}.
\end{equation}

If $\Omega$ is finite, let $p= \max \Omega$ and define $x_p, x_p'$ as before.
Then either $x_p=x_p'$ or $x_p=1-x_p'$, which implies $T_r(x_p)=T_r(x_p')$ and $T_r(x)=T_r(x')$.

If $\Omega$ is infinite, \eqref{eq:tmp-2} implies that for any $\delta>0$, we can find $ p \in \Omega$ such that $|T_r(x_p) -T_r(x_p')| < r^{-p} <\delta$.
Since $\delta$ is arbitrary, this completes the proof. 
\end{proof}

Lemma~\ref{lem:sim} implies that each local level set is contained in a certain level set, that is, $L_{r,x}^{loc} \subset L_r(T_r(x))$.
By the definition of local level sets, we have the following simple facts. Since their proofs are trivial, we omit
the details.

\begin{lemma}\label{lem:lls}
For any given $x \in (0,1)$, we write $n=\card \{ j \in \Z^+ : D_{r,j}(x) =0 \}$. 
\begin{enumerate}
\item If $n =\infty$, then local level set $L_{r,x}^{loc}$ is infinite.
\item If $n <\infty$ and $x $ is an $r$-based rational, then $L_{r,x}^{loc}$ is finite and $\card (L_{r,x}^{loc})= r^{n}$.
\item If $n <\infty$ and $x $ is not an $r$-based rational, then $\card (L_{r,x}^{loc})= 2 r^{n}$.
\item Write $t= \inf L_{r,x}^{loc}$, then $t \in L_{r,x}^{loc}$. Moreover, $D_{r,j}(t) \geq 0$ for all $j \in \Z^+$.
If $t=0.\vep_1 \vep_2 \cdots$, then $\vep_k = r/2$ for each index $k$ satisfying $D_{r,k}(x_0)=0$.
\end{enumerate}

\end{lemma}

\section{Humps}

\subsection{Definition of humps.}
In this subsection, we recall an important concept known as humps, named for their resemblance to the humps on camels' back.

Given $x_0 \in \B^r$,
there exists a unique $m \in \Z^+$ such that $r^{2m} x_0 \in \Z$ and $D_{r,2m}(x_0)=0$. 
This uniqueness follows from the definition of $\B^r$.
Define
$$
 H(x_0) := \gtr \cap \big(I(x_0) \times \R\big),
$$
where $I(x_0) = [x_0,x_0+ r^{-2m}]$.
We call $H(x_0)$ a \emph{hump} of \emph{order} $m$. 
By \eqref{eq:simi-eq}, we know that $H(x_0)$ is a similar copy of the full graph $\gtr$ .
More precisely, there exists a similitude $\Psi$ such that $H(x_0) = \Psi(\gtr)$, where 
\begin{equation}\label{eq:simi}
\Psi(x,y)=\big(x_0+x/r^{2m},T_r(x_0) +y/ r^{2m}\big), \quad (x,y) \in [0,1] \times \R.
\end{equation}
The \emph{generation} of the hump $H(x_0)$ is given by 
$$\card \big\{ j \in \Z^+ : D_{r,j}(x_0) = 0 \big\}.$$
Moreover, if $D_{r,j}(x_0) \geq 0$ for every $1 \leq j \leq 2m$, and $\vep_k = r/2$ for each index $k$ satisfying $D_{r,k}(x_0)=0$, then we say $H(x_0)$ is a \emph{leading hump}.
By convention, we call $\gtr$ a leading hump of generation $0$ and order $0$.

We denote the set of all humps by 
$\H^r := \{ H(x) : x \in \B^r \},$ 
and by $\hath^r$ the subset of $\H^r$ consisting of all leading humps. 
Given $m, n \geq 0$, we denote the set of leading humps of order $m$ and generation $n$ by $\hath_{m,n}^r$. 
Additionally, we have $\hath_{m,0}^r = \hath_{0,m}^r = \emptyset$ for all $m \in \Z^+$, but $\hath_{0,0}^r=\{ \gtr\}$.

\begin{figure}
    \centering
    \begin{minipage}[t]{0.48\linewidth}
           \includegraphics[width=\linewidth]{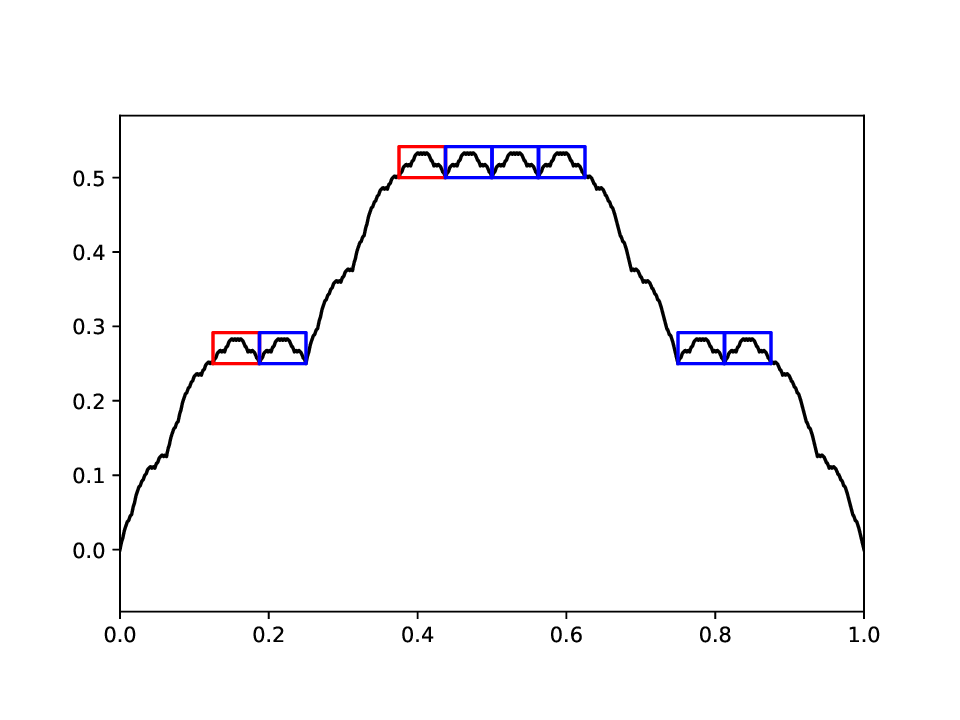}
    \end{minipage}
    \begin{minipage}[t]{0.48\linewidth}
           \includegraphics[width=\linewidth]{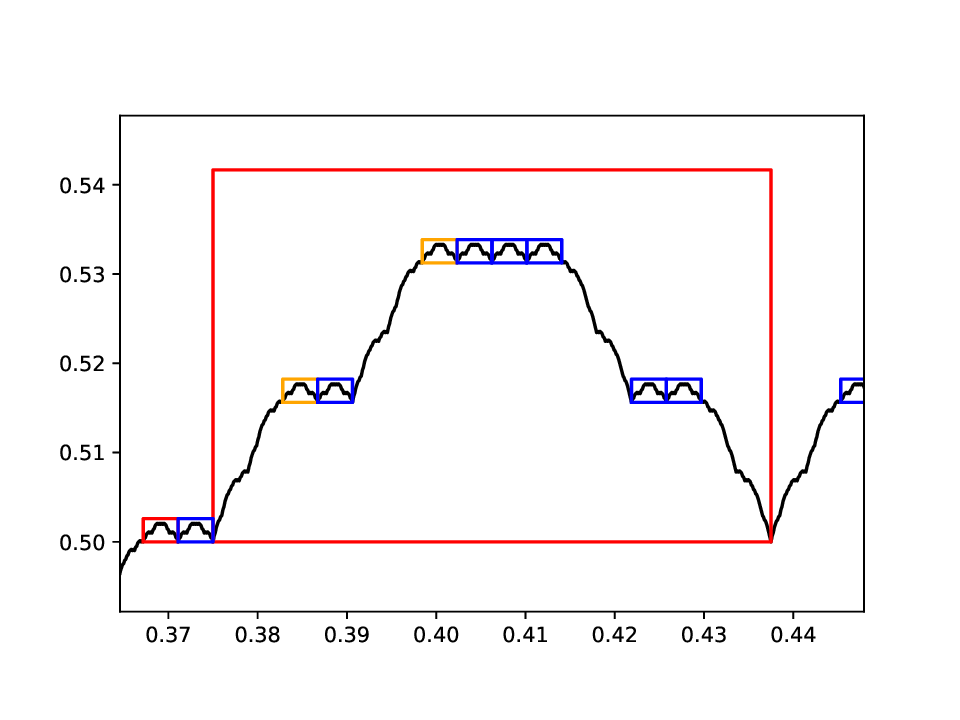}
    \end{minipage}
    \caption{
    Let $r=4$.
    The left picture shows all the humps of order $1$. The right picture displays the humps close to the hump $H(3/16)$.
    The red humps are leading humps with generation $1$, while the orange ones are leading humps with generation $2$.
    }\label{fig:humps}
\end{figure}

\begin{remark}
See Figure \ref{fig:humps} for an illustration of these concepts. 
Assume that $x \in \B^r$ and $H(x)$ is a leading hump.
Consider the set $\{ H(t) : t \in L_{r,x}^{loc} \}$, all the humps contained share the same order and generation.
By the proof of Lemma~\ref{lem:sim}, we have $x=\min L_{r,x}^{loc}$ so that $H(x)$ is the leftmost hump among them.
\end{remark}

Let $x_0 \in \B^r$ and $\Psi$ be as in \eqref{eq:simi}, then the \emph{truncated hump} at $x_0$ is defined by
$$H^t(x_0) = \Psi(\gtrss).$$
Define the projection $\pi_Y$ by letting $\pi_Y(x,y)=y$ for all $(x,y) \in \R^2$.
Now we establish the connection between $N_r^{loc}(y)$ and truncated humps as follows.

\begin{lemma}\label{lem:card}
Given $y >0$, if $L_r(y)$ is finite, then
$$
N_r^{loc}(y) =\card \{ H \in \hath^r : y \in \pi_Y(H^t) \}.
$$
\end{lemma}
\begin{proof}
It is clear that there are finitely many local level sets contained within $L_r(y)$.

Assume $H(x_0)$ is a leading hump of order $m$ such that $y \in \pi_Y(H^t(x_0))$, 
set $ \Lambda = \{ x \in (0,1) : (x,y) \in H^t(x_0)  \} $.
Combining this with Lemma~\ref{lem:tmp} and \eqref{eq:simi},
there are two elements $ \Lambda =\{ x' ,x''\} $, where $x'+x''=2x_0+r^{-2m}$. 
This implies that $x' \sim_r x''$ and $\Lambda \subset L_{r,x'}^{loc}$.
Hence, $\card \{ H \in \hath^r : y \in \pi_Y(H^t) \} \leq N_r^{loc}(y) $.

Conversely, given $x \in [0,1]$ satisfying $L_{r,x}^{loc} \subset L_r(y)$, write $\Omega_x = \{ j : D_{r,j}(x) = 0 \}$.
Set $t =\inf L_{r,x}^{loc}$. Then by Lemma~\ref{lem:lls}, $D_{r,j}(t) \geq 0$ for all $j \in \Z^+$ and $\vep_j=r/2$ for any $j \in \Omega_x$.
From Lemma~\ref{lem:lls} again, $\Omega_x$ is finite. 
We can therefore take $m =\max \Omega_x $.
Define $x_0 = 0.\vep_1\cdots\vep_{m}$ as the $m$-digit truncation of the expansion of $t$ so that $H(x_0) $ is a leading hump.
Since $D_{r,j}(t) >0$ for all $j>m$, we have $(x,y) \in H^t(x_0)$.
As established earlier, $H^t(x_0)$ intersects only one local level set within $L_r(y)$.
This implies $N_r^{loc}(y) \leq \card \{ H \in \hath^r : y \in \pi_Y(H^t) \}$. 
Now the proof is complete.
\end{proof}

\begin{lemma}\label{lem:leb}

Assume that $H(x)$ is a hump of order $m$.
Then
$$\leb \big( \pi_Y(H^t(x))\big)  = \frac{1}{2r^{2m}},$$
where $\leb(\cdot)$ represents the Lebesgue measure on $\R$. 

\end{lemma}
\begin{proof}

Let $S_{r}^*$ be the function defined in the proof of Lemma~\ref{lem:tmp}, then $S_r^*$ is continuous and increasing on $[0,1/2]$.
We obtain $S_r^*(0) = 0 $ and $S_r^*(1/2) = 1/2 $. 
Combining with the symmetry of $\gtr$, we have $\pi_Y (\g S_r^*)= [0,1/2]$ and
$$
	\pi_Y (\g S_r^*) = \pi_Y \big(\g S_r^* |_{X_r^*} \big) \cup \pi_Y \big(\g S_r^* |_{B_1^r} \big)
	= \pi_Y (\gtrss) \cup  T_r (\B_1^r).
$$
Moreover, $T_r(t) \in \mathbb{Q}$ for any $t \in \B_1^r $.
Since $\leb(\mathbb{Q})=0$, it follows that $\leb(T_r(\B_1^r))=0$ and $\leb (\pi_Y (\gtrss) )=1/2$.

Fix $x_0 \in \B^r$ with order $m$, let $\Psi$ be the similitude mentioned in \eqref{eq:simi}.
Thus, we have
$
	\leb \big( \pi_Y(\gtrss) \big) = r^{2m} \leb \big( \pi_Y(H^t(x_0)) \big) .
$
Then the proof is complete.
\end{proof}

\subsection{Cardinality of humps.}

In this subsection, we aim to calculate the cardinality of $\hath_{m,n}^r$.
Firstly, we focus on the case where $r=2$.
Recall that the Catalan number $C_n$ is defined as $C_n = \frac{1}{n+1} \binom{2n}{n}$ for all $n \geq 0$. See \cite[Chapter 2.6.6]{HHM08} for more details.

\begin{lemma}\label{lem:card-n-1}
For any $m \in \Z^+$ we have $\card(\hath_{m,1}^2) = C_{m-1}$.
\end{lemma}
\begin{proof}
Given a binary rational $x=0.\vep_1 \cdots \vep_{2m}$ with $D_{2,2m}(x)=0$, then for any $1 \leq j \leq 2m$, we have
$$
	D_{2,j}(x)=\card\{ 1 \leq k \leq j  : \vep_k=0 \}-\card\{ 1 \leq k \leq j  : \vep_k=1 \}.
$$
Hence, we can regard the cardinality of $\hath_{m,1}^2$ as the number of possible random walks \cite{F71} on $\Z$,
by setting the start point $s_0=0$, with steps $+1$ or $-1$ and 
$$
	s_j :=D_{2,j}(x) = \sum_{k=1}^j  (-1)^{\vep_k} = s_{j-1} +(-1)^{\vep_j}.
$$
Thus, $H(x)$ is a leading hump with generation one if and only if $s_j>0$ for all $1 \leq j < 2m$ and $s_{2m}=0$.
It is clear that $\vep_1=0$ and $\vep_{2m}=1$.
Thus, it is classical that $\card(\hath_{m,1}^2)=C_{m-1}$ \cite{HHM08}. Then the proof is complete.
\end{proof}

\begin{lemma}
For any integers $ m \geq n \geq 2 $, we have
\begin{equation}\label{eq:gen-hath}
	\card(\hath_{m,n}^2)=\sum_{k=1}^{m-n+1} \card(\hath_{k,1}^2) \card(\hath_{m-k,n-1}^2).
\end{equation}
\end{lemma}
\begin{proof}
For any $x \in \B^2$ such that $H(x) \in \hath_{m,n}^2$, assume that its binary expansion is given by $ x=0. \vep_1 \cdots \vep_{2m} $.
Set 
$$p =\min \{1 \leq  j \leq 2m : D_{2,j}(x)=0 \} /2 ,$$ then $p \in \Z^+$.
Write $x_1=0.\vep_1 \cdots \vep_{2p}$ and $x_2=0.\vep_{2p+1} \cdots \vep_{2m}$,
then $H(x_1) \in \hath_{p,1}^2$ and $H(x_2) \in \hath_{m-p,n-1}^2$.
Since $x$ is arbitrary, we get the desired result.
\end{proof}

The following result plays an important role in the proof of Theorems~\ref{thm-1} and~\ref{thm-2}.
\begin{lemma}\label{lem:gen-1-4}
For any $n \in \Z^+$, we have
\begin{equation}\label{eq:s-n-cal}
	\sum_{m=n}^\infty \frac{\card(\hath_{m,n}^2)}{ 4^m} = \frac{1}{2^n}.
\end{equation}
\end{lemma}

\begin{proof}
Fix $n \in \Z^+$, write the generating function \cite[Chapter 2.6]{HHM08} for $\{ \card(\hath_{m,n}^2) \}$ as
$$
	P_n(z):= \sum_{m=n}^\infty \card(\hath_{m,n}^2) z^m, \quad  z \in \R. 
$$
For any $n \in \Z^+$, we claim that $P_n(z)=(P_1(z))^n$. It suffices to show that $P_n(z)=P_1(z)P_{n-1}(z)$. The coefficients of $z^m$ on both sides coincide by \eqref{eq:gen-hath}, so the two power series are identical, proving the claim.

From Lemma~\ref{lem:card-n-1}, we have $\card(\hath_{m,1}^2)= C_{m-1}$. 
The generating function for $\{C_{m-1}\}$ is known as
$$
P_1(z)  = \sum_{m=1}^{\infty} C_{m-1} z^m = \frac{1 - \sqrt{1 - 4z}}{2}, \quad z \in [0, 1/4].  
$$
Hence $P_n(z) = (\frac{1 - \sqrt{1 - 4z}}{2})^n$. Taking $z=1/4$
 yields the result.
 \end{proof}

The following result is classical.
\begin{lemma}[\cite{A12,AK11}]\label{lem:card-tradition}
Let $r=2$ and $m \in \Z^+$. There are $\binom{2m}{m}$ humps of order $m$. 
Among these, there are $C_m$ leading humps of order $m$.
\end{lemma}

\begin{lemma}\label{lem:count-old}
For any $m \in \Z^+$ and even integer $r$, there are $ (\frac{r}{2})^{2m} \binom{2m}{m} $ humps of order $ m $.
\end{lemma}

\begin{proof}
For any binary balanced rational $x= 0. \vep_1 \cdots  \cdots \vep_{2m}$, we can transform it into an $r$-based balanced rational in the following way. Every ``$0$'' in the binary expansion can be replaced by one digit in $\{0,1,\ldots, r/2-1\}$, while every ``$1$" can be turned into $\{ r/2,\ldots, r-1\}$. Obviously this establishes an $1$-to-$(r/2) ^{2m}$ correspondence between the binary humps and $r$-based humps with order $m$.
\end{proof}

Applying a similar method, we can count the leading humps.

\begin{lemma}\label{lem:count-tmp}
For any $m ,n\in \Z^+$ and even integer $r$, we have
\begin{equation}
\card(\hath_{m,n}^r) = \Big( \frac{r}{2} \Big)^{2m-n} \card(\hath_{m,n}^2).
\end{equation}
\end{lemma}

\begin{proof}

By the definition of leading humps, if $D_{r,j}(x)=0$ then $\vep_j$ must be $r/2$. As for a leading hump of generation $n$, there are $n$ such digits without the freedom of choice. Only the rest $2m-n$ digits can each have $r/2$ choices.
\end{proof}

\section{Proof of the main result}

In this section, we write the maximum value of $T_r$ as $M_r$.
Baba \cite{B84} proved that $M_r=r^2 / (2r^2 - 2)$ for any even integer $r\geq 2$.
Now we can prove the main result.

\begin{proof}[Proof of Theorem~\ref{thm-1}]
Allaart \cite{A14} showed that $L_r({t})$ is finite for almost every $t$.
Since $y$ is uniformly distributed on $[0, M_r]$, 
by the definition of expectation and Lemma~\ref{lem:card},
$$
	\E \big[N_r^{loc}(y) \big] = \frac{1}{M_r}\int_0^{M_r} N_r^{loc}(y) \, dy 
	=  \frac{1}{M_r} \int_0^{M_r} \card  \big\{ H \in \hath^r : y \in \pi_Y(H^t) \big\} \, dy.
$$
Moreover,
$$
	 \int_0^{M_r} \card \big\{ H \in \hath^r : y \in \pi_Y(H^t) \big\} \, dy
	 = \sum_{H \in \hath^r} \leb \big( \pi_Y(H^t) \big)  .
$$

Note that $\hath^r$ can be written as the disjoint union of the following sets
$$
	 \hath^r 	= \bigcup_{n=0}^\infty \bigcup_{m=n}^\infty \hath_{m,n}^r.
$$
For $n =0$, recall that $\gtr$ is the only hump of generation $0$. 
Then by Lemma~\ref{lem:leb}, $\leb ( \pi_Y(H^t(0)) )=\leb ( \pi_Y(\gtrss) )=1/2 $.
By Lemmas~\ref{lem:gen-1-4} and \ref{lem:count-tmp}, for any $n \in \Z^+$,
$$
	 \sum_{m=n}^\infty \sum_{H \in \hath_{m,n}^r}  \leb \big(  \pi_Y(H^t) \big) 
	  =\sum_{m=n}^\infty  \frac{\card(\hath_{m,n}^r) }{2r^{2m}} 
	=\sum_{m=n}^\infty   \frac{2^n \card(\hath_{m,n}^2)}{2r^n \cdot 4^m} 
	= \frac{1}{2 r^n} .
$$
Hence,
$$
	 \sum_{H \in \hath^r} \leb \big( \pi_Y(H^t) \big) 
	= \sum_{n=0}^\infty \sum_{m=n}^\infty \sum_{H \in \hath_{m,n}^r} \leb \big( \pi_Y(H^t) \big) 
	=  \sum_{n=0}^{\infty} \frac{1}{2 r^n} =\frac{r}{2(r-1)}.
$$
Therefore, we have
$$
	\E \big[N_r^{loc}(y) \big] 
	=  \frac{1}{M_r} \sum_{H \in \hath^r} \leb \big( \pi_Y(H^t) \big) 
	= \frac{r}{2(r-1) M_r} =1+ \frac{1}{r} .
$$
Thus, the proof is complete.
\end{proof}

It is natural to ask whether the expected cardinality of a level set $L_r(y)$ is also finite.
The following theorem provides a negative answer.

\begin{theorem}\label{thm-2}
Let $y$ be a random variable uniformly distributed on $[0, M_r]$. Then the expected cardinality of its level set $ L_r(y) $ is infinite. That is,
$$
	\E \big[\card \big( L_r(y)\big) \big] = \frac{2r^2-2}{r^2}  \int_0^{\frac{r^2}{2r^2-2}} \card  \big( L_r(y)\big) \, dy = \infty.
$$
\end{theorem}
\begin{proof}
Using an argument similar to that in Theorem~\ref{thm-1}, we have
$$
	\E\big[ \card\big(L_r(y)\big)\big]
	= 2 \E \big[ \card \big\{ H \in \H : y \in \pi_Y(H^t) \big\} \big] 
	= 2 \sum_{H \in \H} \leb\big(   \pi_Y(H^t) \big) .
$$
Combining this with Lemmas \ref{lem:leb} and \ref{lem:count-old} yields
$$
	\E\big[ \card\big(L_r(y)\big)\big]
	= 2 \sum_{m=0}^{\infty} \Big( \frac{r}{2} \Big)^{2m}  \binom{2m}{m}  \cdot \frac{M_r}{r^{2m}}  
	=2 M_r  \sum_{m=0}^{\infty}  \frac{1}{4^m }      \binom{2m}{m}    .
$$
Since $\binom{2m}{m} \sim 4^m / \sqrt{  m}$, it follows that $\E\big[ \card\big(L_r(y)\big)\big]= \infty$.
\end{proof}

\subsection*{Acknowledgements.}  

The authors would like to thank Prof. Shi-Lei Kong, Prof. Yan-Qi Qiu and Prof. Huo-Jun Ruan for helpful discussions. 
Ying and Zhang were financially supported by Zhejiang University Student Research Training Program No. 2024R401135.

\end{document}